\documentclass[letterpaper, 10 pt, conference]{ieeeconf}  

\IEEEoverridecommandlockouts                              

\usepackage{graphics} 
\usepackage{amsmath} 
\usepackage{amssymb}  
\usepackage{graphicx}
\usepackage{algorithm}
\usepackage{algorithmic}
\usepackage{comment}

\usepackage[margin=0.75in]{geometry}

\usepackage{algorithm}
\usepackage{xcolor}
\usepackage{comment}

\DeclareRobustCommand{\bigO}{%
  \text{\usefont{OMS}{cmsy}{m}{n}O}%
}

\newcommand{\revision}[1]{\textcolor{black}{#1}}

\allowdisplaybreaks

\title{\LARGE \bf
Distributed Algorithm for Linear Least Squares Problems via Iteratively Preconditioning by Successive Over-Relaxation
}

\author{Tianchen Liu, \IEEEmembership{Member, IEEE}, Kushal Chakrabarti, \IEEEmembership{Member, IEEE}  
\thanks{Tianchen Liu is with State Key Laboratory of Robotics and Systems (HIT), Harbin, 150000 China. (Email: tianchenliu@hit.edu.cn)} 
\thanks{Kushal Chakrabarti is with the Division of Data \& Decision Sciences, Tata Consultancy Services Research, Mumbai, 400607 India. (Email: chakrabarti.k1@tcs.com)} 
\thanks{This work is supported by Self-Planned Task (No. SKLRS202605C) of State Key Laboratory of Robotics and Systems (HIT).} 
\thanks{Corresponding Author: Tianchen Liu.} 
}

\begin{document}

\maketitle
\thispagestyle{empty}
\pagestyle{empty}

\begin{abstract}
\revision{This paper focuses on a novel distributed optimization algorithm for solving multi-agent linear least squares problems within a server-agent framework. Distributed methods are particularly suited to large-scale systems. The computational tasks can be partitioned across agents, while local data of one agent remain inaccessible to others. The proposed approach employs an iterative preconditioning scheme with the successive over-relaxation technique, leading to a substantial reduction of the number of iterations required for convergence, which presents a critical advantage where server-agent communication is costly.} Theoretical convergence proof of the proposed algorithm is presented. Finally, benchmark problems are tested to demonstrate the performance comparison of the proposed approach over several prominent distributed methods.
\end{abstract}

\section{Introduction}
\label{sec:introduction}
Distributed optimization approaches can help in effectively addressing large-scale interconnected problems by dividing complex problems into multiple computing nodes. These algorithms leverage the computational power of multiple machines or agents, leading to more efficient use of available resources. This parallel processing capability allows for tackling problems that are computationally intensive to solve in a centralized manner, and is also beneficial in engineering applications where local computations reduce the need for extensive data transmission (e.g., sensor networks, multi-robot systems, etc.). In addition, in scenarios where data are distributed across multiple locations or owned by different entities, distributed methods enable collaborative problem-solving, while local data of one entity remain inaccessible to others. This property is critical in preserving local sensitive information against unauthorized access.

In this paper, an accelerated distributed optimization algorithm is proposed to solve an $N$-agent linear least squares (LS) problem within a server-agent structure. Suppose that all agents share the same \textit{parameter vector} $x \in \mathbb{R}^d$. For each agent $i=1, \dots, N$, $A_i \in \mathbb{R}^{n_i \times d}$ is a real-valued \textit{local data matrix} representing a set of $n_i$ data points, and $b_i \in \mathbb{R}^{n_i}$ is a real-valued \textit{local observation vector}. A \textit{local cost function} $F_i(x) = \frac{1}{2} \| A_i x - b_i \|^2$ is defined to update $x$ locally for agent $i$ using its local data. 

\begin{figure}[ht!]
    \centering
    \includegraphics[width=0.7\linewidth]{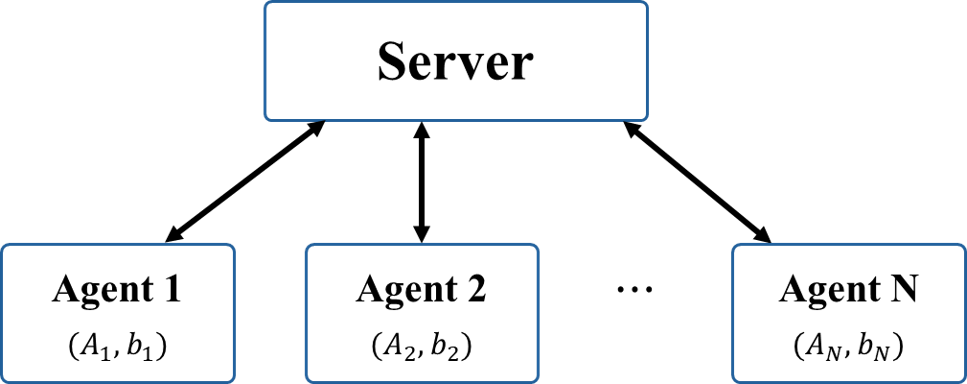}
    \caption{Distributed system structure.}
    \label{fig:dist_sys}
\end{figure}



Since each agent cannot access the data of others, a central server facilitates bidirectional communication with all agents, ensuring the coordination of updates among them. The overall structure is referred to as a {\em server-agent} framework, as shown in Fig.~\ref{fig:dist_sys}. Hence, we formulate the overall problem by defining the \textit{aggregate cost function} as
\begin{align}
\resizebox{0.36\textwidth}{!}{$
F(x) = \sum_{i=1}^{N} F_i(x) = \sum_{i=1}^{N} \frac{1}{2} \| A_i x - b_i \|^2,
$}
\label{eq:overall_obj_fun}
\end{align}
which is the sum of all local cost functions for agents $i=1, \dots, N$. The goal is to compute the optimal parameter vector $x^* \in \mathbb{R}^d$ such that $x^* = \arg \min_{x} F(x)$. The terminology \textit{distributed} refers to the process where updates are distributed across agents \textcolor{black}{and updates are computed locally by each agent using its dataset $(A_i,b_i)$ and shared solution $x$.}

The multi-agent linear LS problem presented in~\eqref{eq:overall_obj_fun} has a wide range of applications, such as data fitting~\cite{zhang2017stochastic}, state estimation~\cite{sasai2023linear}, and LiDAR simultaneous localization and mapping~\cite{pan2021mulls}. Moreover, regarding the machine learning field, various problems can be framed as linear LS, including chain demand forecasting~\cite{carbonneau2008application}, prediction of online user inputs~\cite{canny2013method}, etc. 
\revision{Specifically, one motivation is originated from the federated learning concept~\cite{smith2017federated}, where data are processed locally by each source/agent, and inaccessible to others.} This leads to a necessity of developing distributed algorithms for solving data-driven modeling problems in a server-agent framework.
\revision{If fast convergence can be achieved in such applications (e.g., federated learning or sensor networks), communication rounds and energy consumption can be directly reduced, enabling more efficient model training and problem solving.}

\textcolor{black}{In many distributed learning environments, the communication cost, rather than computation at the local level, could be the bottleneck. Each iteration may require a communication round across slow, unreliable, or high-latency links. Even if one round involves heavy local computation, it’s still better than requiring hundreds of synchronization rounds. Such environments include federated learning across mobile users (slow global syncing), cross-continent or inter-data-center machine learning training (high-latency links with round trips cost ranging from milliseconds to seconds), and federated medical data networks (hospitals cannot transmit frequently). Similarly, power system operators run real-time distributed optimization (e.g., optimal power flow, demand-response coordination) over geographically spread controllers. Iteration count matters because communication delays between grid substations are high and real-time optimization must converge within strict deadlines (e.g., 60 Hz control loops, emergency response). Hence, algorithms like ADMM with expensive per-iteration local complexity are preferred because they converge in few iterations over lightweight algorithms like vanilla gradient-descent. Robots or drones collaborating on mapping or coordination tasks often use heavy algorithms like second-order consensus or distributed Newton methods to reduce number of global rounds, due to battery limitation, need for fast consensus, and intermittent communication. Privacy-preserving distributed optimization, using homomorphic encryption or secure multiparty computation, has already high computational cost per iteration due to cryptographic operations, but their communication overhead is even more expensive. Finally, distributed controllers in automotive systems with MPC, control of autonomous vehicles in traffic platoons, or robotics manipulation across multiple controllers may only have limited time for a fixed iterations per control cycle. In all the aforementioned problems, communication or synchronization is the bottleneck in terms of latency, energy, privacy, or bandwidth. Thus, it is preferable to use a distributed optimization algorithm with expensive local iterations but low iteration count to solve such problems. In this paper, we propose a distributed algorithm for solving~\eqref{eq:overall_obj_fun} ideally suited to such scenarios.}

Our proposed algorithm is an extension of our previous works~\cite{chakrabarti2022iterative, liu2024IpgSor}. 
In~\cite{chakrabarti2022iterative}, the {\em IPG} (Iteratively Preconditioned Gradient-descent) algorithm is introduced to solve~\eqref{eq:overall_obj_fun} in a server-agent framework. \textcolor{black}{IPG is a second-order method guaranteed to converge in fewer iterations than gradient-descent and its momentum variants~\cite{chakrabarti2022iterative}.} Our aim is to propose a novel accelerated distributed algorithm to improve the performance for~\eqref{eq:overall_obj_fun} \textcolor{black}{primarily in terms of reduced number of iterations.} Following this idea, we propose a distributed version of the {\em IPG-SOR} ({\em IPG} using successive over-relaxation (SOR)) method~\cite{liu2024IpgSor} in this paper.
\revision{Compared with the distributed {\em IPG} method~\cite{chakrabarti2022iterative}, the incorporation of the SOR technique removes step size computation at each iteration, and can significantly reduce the iteration counts for convergence. It is worth noting that although this work is the distributed variant of our centralized {\em IPG-SOR} method~\cite{liu2024IpgSor}, its convergence conditions and analytical proofs differ substantially to accommodate the distributed structure.} In summary, the main contributions are as follows:
\begin{itemize}
\item \revision{A novel accelerated distributed algorithm is proposed to solve multi-agent linear LS problems, significantly reducing the number of iterations needed for convergence.}
\item Theoretical convergence proof of the algorithm is presented.  
\item Numerical results of benchmark problems are tested and analyzed to validate the efficacy of the proposed algorithm. 
\end{itemize}

The rest of the paper is organized as follows: Section II introduces the background of {\em IPG}, {\em IPG-SOR} and related work. The proposed approach and the theoretical results are presented in Section III and Section IV, respectively. In Section V, numerical examples are tested to validate the algorithm. Section VI concludes the paper.

\section{Background}

\subsection{IPG \& IPG-SOR}
A brief introduction of the {\em IPG} algorithm (distributed~\cite{chakrabartiAcceleratingDistributedConvex2021} or centralized~\cite{liuIterativelyPreconditionedGradientDescent2023}) proposed in our previous work is as follows. Let $x \in \mathbb{R}^{d}$ be the optimizing variable vector to minimize the convex objective function $F(x)$. A preconditioner matrix $K \in \mathbb{R}^{d \times d}$ is introduced to help the iterative update of $x$ to converge to the optimum $x^{*}$. Before the iterations start, initial values of $x(0)$, $K(0)$, and two scalar constants $\beta \geq 0$, $\delta > 0$ are selected. 
At the $t^{\text{th}}$ iteration ($t = 0, 1, \ldots$), $x(t)$ and $K(t)$ are updated as
\begin{align}
x(t+1) & = x(t) - \delta K(t) g(t), \label{eq:update_x_ipg}\\
K(t+1) & = K(t) - \alpha(t) [H(t) K(t) + \beta K(t) - I],\label{eq:update_K_ipg}
\end{align}
until $\| x(t+1) - x(t) \|< \varepsilon$, where $\varepsilon$ is a predetermined tolerance value. Here, $g(t) \in \mathbb{R}^{d}$ and $H(t) \in \mathbb{R}^{d \times d}$ denote the gradient and Hessian of $F(x)$ with respect to $x$, evaluated at $x = x(t)$. $I \in \mathbb{R}^{d \times d}$ denotes the identity matrix. $\alpha(t)$ is the step size,  calculated by convergence conditions (see~\cite{chakrabartiAcceleratingDistributedConvex2021}). 

The {\em IPG-SOR} method~\cite{liu2024IpgSor} is derived based on {\em IPG}, with the preconditioner $K$ updated by the SOR technique~\cite{kelleyIterativeMethodsOptimization1999}. Two scalars, $\omega \in (0, 2)$ and $\beta \geq 0$, are selected. At iteration $t$, let $P_{\omega}(t) = \frac{1}{\omega} (D(t) + \beta I) + L(t) \in \mathbb{R}^{d \times d}$, where $D(t) \in \mathbb{R}^{d \times d}$ and $L(t) \in \mathbb{R}^{d \times d}$ are square matrices containing the diagonal and the strictly lower parts of $H(t)$, respectively. Given initial $x(0), K(0)$, $x(t)$ and $K(t)$ are updated as
\begin{align}
x(t+1) & = x(t) -  K(t) g(t),\label{eq:update_x_sor} \\
K(t+1) & = K(t) - P_{\omega}^{-1}(t) [H(t) K(t) + \beta K(t) - I]. \label{eq:update_K_sor}
\end{align} 
The iterative process stops when $\| x(t+1) - x(t) \|< \varepsilon$. 

Both {\em IPG} and {\em IPG-SOR} can be regarded as quasi-Newton methods. Their faster convergence performance than several other approaches (including Nesterov accelerated gradient-descent (NAG)~\cite{nesterov1983method}, Adam~\cite{kingma2014adam}, BFGS~\cite{broydenConvergenceClassDoublerank1970}, etc.) for solving unconstrained convex problems have been presented in~\cite{chakrabartiAcceleratingDistributedConvex2021, liuIterativelyPreconditionedGradientDescent2023}.

\subsection{Related Work}
For the distributed optimization problem over the server-agent framework, multiple approaches have been developed. The literature such as~\cite{lessard2016analysis} has discussed the deployment of the concept of momentum from the NAG algorithm~\cite{nesterov1983method} to significantly accelerate the original gradient-descent method in a distributed manner. The APC method~\cite{azizan2019distributed} is a combination of the NAG method with a projection operation, which shows a faster convergence. 
However, these mentioned algorithms require either $A_i x = b_i\, (i = 1, \ldots, N)$ to have a unique solution, or the solution of~\eqref{eq:overall_obj_fun} to be unique. It is also known that the convergence speed is fundamentally limited by the condition number of the collective data. \revision{Among advanced distributed optimization approaches, the Alternating Direction Method of Multipliers (ADMM)~\cite{boyd2011distributed} is widely used. For linear LS, the overall problem is decomposed into local subproblems with consensus constraints. The distributed L-BFGS method~\cite{najafabadi2017large} is another commonly used one due to the low computational costs, which extends the quasi-Newton techniques in a distributed manner and achieves faster convergence than first-order methods for large-scale problems.}



\section{Proposed Approach}
\label{sec:approach}
In this section, we present a novel distributed algorithm for the multi-agent linear LS problem~\eqref{eq:overall_obj_fun}. The algorithm is derived to solve~\eqref{eq:update_x_sor}, \eqref{eq:update_K_sor} in {\em IPG-SOR} in a distributed manner.

Recall that each agent $i$ has its local data matrix $A_i$ and observation vector $b_i$. Since the server has access to data from all agents, it is assumed that the collective data matrix $\mathbf{A} = [A_1^T \, A_2^T \, \dots \, A_N^T]^T$, which concatenates all $A_i$'s in a vertical sequence, can be obtained. Let $D$ and $L$ be the square matrices containing the diagonal and the strictly lower parts of $\mathbf{A}^T \mathbf{A}$, respectively. Hence, we calculate the matrix
\begin{align}
\resizebox{0.3\textwidth}{!}{$
    P = P_{\omega}^{-1} = \left( \frac{1}{\omega} (D + \beta I) + L \right)^{-1}
    $}
    \label{eq:calc_P_omega_inv}
\end{align}
in the server. Note that if $\mathbf{A}$ remains constant, then $P$ is computed only once at the first iteration.

\begin{algorithm}[ht!]
\caption{\textit{Distributed IPG-SOR} algorithm for multi-agent linear LS problems in a server-agent framework.}
\label{alg:1}
\begin{algorithmic}[1]
    \STATE Server initializes $x(0) \in \mathbb{R}^d$, $K(0) \in \mathbb{R}^{d \times d}$, $\beta \geq 0$, and $\omega \in (0, 2)$.
    \STATE Server obtains collective data $\mathbf{A}$ and computes $P$ by \eqref{eq:calc_P_omega_inv}.
    \FOR{$t = 0, 1, 2, \dots$}
        \STATE Server sends $x(t)$ and $K(t)$ to each agent $i \in \{1, \dots, N\}$.
        \STATE Each agent $i$ computes local gradient $g_i(t)$ by \eqref{eq:calc_g_i}, and a set of vectors $\{ R_{i,j}(t), \, j = 1, \dots, d \}$ by \eqref{eq:calc_r_ij}.
        \STATE Each agent $i$ sends $g_i(t)$ and the set $\{ R_{i,j}(t), \, j = 1, \dots, d \}$ to server.
        \STATE Server updates $K(t)$ to $K(t+1)$ by \eqref{eq:update_k_j}.
        \STATE Server updates $x(t)$ to $x(t+1)$ by \eqref{eq:update_x_ipg_sor_dist}.
    \ENDFOR
\end{algorithmic}
\end{algorithm}

The proposed algorithm is summarized in Algorithm~\ref{alg:1}. Before initiating the iterative process, the server selects initial values $x(0)$ and $K(0)$, as well as two scalar parameters $\beta$ and $\omega$. The parameter $\beta$ is broadcast to all agents. Next, for iteration $t=0, 1, \dots$, the solution steps are as follows.


\subsection{Steps in Iteration $t$}
\begin{figure}[ht!]
    \centering
    \includegraphics[width=0.9\linewidth]{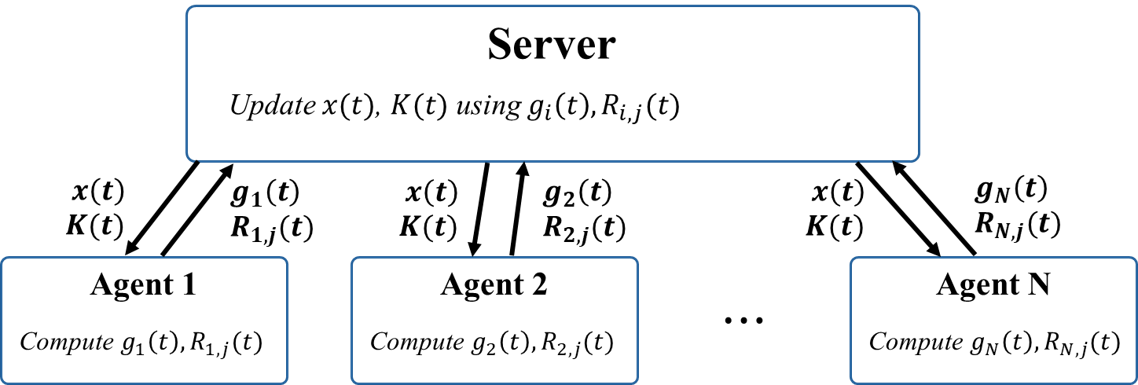}
    \caption{Distributed solution approach at iteration $t$.}
    \label{fig:dist_ipg_sor}
\end{figure}

In the $t^{\text{th}}$ iteration, the algorithm consists of four steps, executed collaboratively by the server and the agents. A graphical illustration is shown in Fig.~\ref{fig:dist_ipg_sor}. 

\textbf{Step 1:} The server sends the current values of $x(t)$ and $K(t)$ to each agent.

\textbf{Step 2:} For each agent $i$, the local gradient $g_i(t)$ of its local
cost function $F_i(x)$ at $x(t)$ is computed as
\begin{align}
    g_{i}(t) = A_i^T \left( A_i x(t) - b_i \right).
    \label{eq:calc_g_i}
\end{align}
Additionally, for each $j = 1, \dots, d$, a set of vectors $\{R_{i,j}(t)\}$ is calculated as follows, 
\begin{align}
\resizebox{0.36\textwidth}{!}{$
    R_{i,j}(t) = \left( A_i^T A_i + \left(\frac{\beta}{N}\right)I \right) k_j(t) - \left(\frac{1}{N}\right)e_j.
    $}
    \label{eq:calc_r_ij}
\end{align}
where $k_j(t)$ and $e_j$ denote the $j^{\text{th}}$ columns of matrices $K(t)$ and $I$, respectively. 

\textbf{Step 3:} Each agent $i$ sends $g_i(t)$ and $\{R_{i,j}(t)\}$ to the server.

\textbf{Step 4:} The server updates the $j^{\text{th}}$ column of the preconditioner matrix $K(t)$ by 
\begin{align}
\resizebox{0.3\textwidth}{!}{$
    k_j(t + 1) = k_j(t) - P \sum_{i=1}^N R_{i,j}(t)
    $}
    \label{eq:update_k_j}
\end{align}
for all $j = 1, \dots, d$, leading to an updated $K(t+1)$. Finally, the parameter vector $x(t + 1)$ is computed by 
\begin{align}
\resizebox{0.3\textwidth}{!}{$
    x(t + 1) = x(t) - K(t + 1) \sum_{i=1}^N g_i(t).
    $}
    \label{eq:update_x_ipg_sor_dist}
\end{align}
The iterative process repeats from \textbf{Step 1} for the next iteration $t+1$ until $\| x(t + 1) - x(t) \| < \varepsilon$.


\subsection{Algorithm Complexity}
\label{subsec:complexity}
Next the operational cost of the proposed algorithm is analyzed, with regard to both computational and communication complexity, which is evaluated by the total number of floating-point operations (flops) required per iteration. We ignore the additions while counting the total number of flops since multiplication is more costly than addition~\cite{golub_matrix_2013}.

1) Computational complexity: Before the beginning of all iterations, the matrix $P$ is calculated only once in the server in $\bigO(d^3)$ flops. For each agent in iteration $t$, the computation of $g_i(t)$ requires $\bigO(n_i d)$ flops, and that of $\{R_{i,j}(t)\}$ requires $\bigO(n_i d^2)$ flops. Note that $\{R_{i,j}(t)\}$ is independently calculated for each agent, hence the computation can be processed in a parallel manner. In the next steps, the server computes $K(t+1)$ in $\bigO(d^2)$ flops, and $x(t+1)$ in $\bigO(d^2)$ flops. 

2) Communication complexity: For each iteration $t$. The server first sends $x(t)$ and $k_j(t)$ to all agents, and an agent sends $g_i(t)$ and $d$ vectors $\{R_{i,j}(t)\}$ back to the server. Hence, the communication complexity of Algorithm 1 for
each agent is $\bigO(d^2)$ flops per iteration, the same as in our previous work in~\cite{chakrabarti2022iterative}. This is more costly than the traditional distributed GD method, which only requires $\bigO(d)$ flops. However, the proposed algorithm can significantly reduce the number of iterations to achieve convergence.

\subsection{\revision{Remarks}}
\revision{An inherent assumption for the proposed approach is that the communication is consistent and synchronous between the server and all agents. Asynchronous communication or network failure cases, which are commonly encountered in distributed systems, have not been considered in this paper. Intuitively, convergence can still be achieved under asynchronous data exchange, but a rigorous proof is not yet available. In addition, the collective data matrix $\mathbf{A}$, as well as local data matrices $A_i$, can be rank deficient. Hence, the feasibility of data distribution can always be satisfied.}

\revision{As seen in Step 4, dense computations are performed in the central server, so that the computation cost per iteration can be much heavier than that of other methods (such as distributed L-BFGS or ADMM). However, given the reduced number of iterations for convergence, when communication overhead dominates the computational cost at each iteration, the reduction in iteration count becomes a decisive factor in improving overall performance. Hence, we conclude that the proposed algorithm are most appropriate for such scenarios, which was our motivation as stated in Section~\ref{sec:introduction}. Besides, a recent collaborating paper on quantum computing~\cite{srinivasan2024trapped} presents an application where reducing iteration counts is of priority, since additional noise occurs at each iteration. \textcolor{black}{Moreover, this high per-iteration cost in the server is the {\em total} cost of Step 4.} Step 4 can be executed in parallel computing units or workers within the server, which can further reduce overall computation time. \textcolor{black}{Training over multiple workers in a single machine (server) is common and cheap in modern GPU units due to availability of cloud services.}}

\section{Convergence Analysis}
In this section, the theoretical results for the convergence of the proposed algorithm in Section~\ref{sec:approach} are presented. 

First, we introduce notations to be used. $\lambda_1, \dots, \lambda_d$ denote the eigenvalues of the positive semidefinite matrix $\mathbf{A}^T \mathbf{A}$ such that $\lambda_1 \geq \dots \geq \lambda_d \geq 0$. $\rho(\cdot)$ denotes the spectral radius of a matrix. Let $\text{Diag}(v)$ denotes a square matrix with the diagonal entries being the elements of the vector $v$ in the same order. For any square matrix $M$, let $\| M \|$ and $\| M \|_F$ denote the \textit{induced 2-norm} and \textit{Frobenius norm} of the matrix, respectively. The collective observation vector concatenating all $b_i$'s in a vertical stack is denoted as $\mathbf{b} = [b_1^T \, b_2^T \, \dots \, b_N^T]^T$. 

Define $K_{\beta} = \left( \mathbf{A}^T \mathbf{A} + \beta I \right)^{-1}$ and $\widetilde{K}(t) = K(t) - K_{\beta}$, then we denote the $j^{\text{th}}$ column ($j \in \{1, \dots, d\}$) of $K(t)$, $K_\beta$ and $\widetilde{K}(t)$ as $k_{j}$, $k_{j\beta}$ and $\widetilde{k}_j(t)$, respectively. Let $g(t)$ be the gradient of the aggregate cost function $F(x)$ at $x(t)$. For iteration $t$, since $\sum_{i=1}^{N} A_i^T A_i = \mathbf{A}^T \mathbf{A}$ and $\sum_{i=1}^{N} A_i^T b_i = \mathbf{A}^T \mathbf{b}$, $g(t)$ by~\eqref{eq:calc_g_i} can be computed as
{\small
\begin{align}
    g(t) & = \textstyle\sum_{i=1}^{N} g_i(t) = \left( \textstyle\sum_{i=1}^{N} A_i^T A_i \right) x(t) - \left( \textstyle\sum_{i=1}^{N} A_i^T b_i \right) \nonumber\\
    & = \mathbf{A}^T (\mathbf{A} x(t) - \mathbf{b}).
\end{align}
}
Note that the solution for the problem~\eqref{eq:overall_obj_fun} can be non-unique. \textcolor{black}{Specifically, the solution is determined by the total gradient $g(t) = \mathbf{0}_d$.} Hence, the least-squares solution set $X^*$ is defined as $X^* = \{ x^* \in \mathbb{R}^d : \textcolor{black}{\mathbf{A}^T}(\mathbf{A} x^* - \mathbf{b}) = \mathbf{0}_d \}$. 

\subsection{Convergence Proof}
First, necessary lemmas for the convergence proof of Algorithm 1 are introduced.

\textbf{Lemma 1.} 
Let $\omega \in (0, 2)$ be selected and $P$ be calculated in~\eqref{eq:calc_P_omega_inv}, then there exists a real parameter $\rho_{\omega} \in (0, 1)$ such that
\begin{align}
\rho(I - P (\mathbf{A}^T \mathbf{A} + \beta I)) \leq \rho_{\omega} < 1.
\end{align}

\textit{Proof:} The updating equation for the preconditioner $K$ (i.e.,~\eqref{eq:update_K_sor}) is
\begin{align}
K(t+1) & = (I - P (\mathbf{A}^T \mathbf{A} + \beta I) ) K(t) - P.
\label{eq:update_K_sor_new}
\end{align}
When $\beta \geq 0$, since $\mathbf{A}^T \mathbf{A} + \beta I$ is symmetric and positive definite, $\omega \in (0,2)$ is the necessary and sufficient condition for $K(t)$ converging to $K_{\beta}$ using~\eqref{eq:update_K_sor_new} (Theorem~3.41 of~\cite{hackbuschIterativeSolutionLarge2016}). Hence,~\eqref{eq:update_K_sor_new} is a contraction, which is equivalent to $\rho(I - P (\mathbf{A}^T \mathbf{A} + \beta I)) < 1$. Hence, a positive parameter $\rho_{\omega} \in (0, 1)$ can be found such that $\rho(I - P (\mathbf{A}^T \mathbf{A} + \beta I)) \leq \rho_{\omega}$. \hfill $\blacksquare$

\textbf{Lemma 2.} Given $\rho_{\omega} \in (0, 1)$ in Lemma 1, then for all $j = 1, \dots, d$, we have
\begin{align}
\|k_j(t+1) - k_{j\beta}\| & \leq \rho_{\omega} \|k_j(t) - k_{j\beta}\|, \quad \forall t \geq 0.
\label{eq:bound_tilde_k_j}
\end{align}

\textit{Proof:} 
Since $\mathbf{A}^T \mathbf{A} = \sum_{i=1}^{N} A_i^T A_i$, by summing up $R_{i,j} (t)$ from $i=1$ to $N$ in~\eqref{eq:calc_r_ij}, we have 
\begin{align}
\textstyle\sum_{i=1}^N R_{i,j} (t) &= \left[(\mathbf{A}^T \mathbf{A} + \beta I) k_j(t) - e_j\right].
\end{align}
Hence, for all $j = 1, \dots, d$,~\eqref{eq:update_k_j} is rewritten as, 
\begin{align}
k_j(t+1) & = k_j(t) - P \left(\mathbf{A}^T \mathbf{A} + \beta I\right) k_j(t) + P e_j.
\label{eq:update_k_j_sum}
\end{align}
Since $K_{\beta} = \left( \mathbf{A}^T \mathbf{A} + \beta I \right)^{-1}$, $\left(\mathbf{A}^T \mathbf{A} + \beta I\right) k_{j\beta} = e_j$. Then 
\eqref{eq:update_k_j_sum} becomes
\begin{align}
\widetilde{k}_j(t+1) 
& = k_j(t) - P \left(\mathbf{A}^T \mathbf{A} + \beta I\right) k_j(t) + P e_j - k_{j\beta} \nonumber\\
& = \widetilde{k}_j(t) - P \left(\mathbf{A}^T \mathbf{A} + \beta I\right) \left(k_j(t) - k_{j\beta}\right) \nonumber\\
& = \left[I - P \left(\mathbf{A}^T \mathbf{A} + \beta I\right)\right] \widetilde{k}_j(t).
\end{align}
From {\it Lemma 1}, $\|I - P(\mathbf{A}^T \mathbf{A} + \beta I)\| \leq \rho_{\omega}$, hence
$\| \widetilde{k}_j(t+1) \| \leq \rho_{\omega} \|\widetilde{k}_j(t)\|$, which indicates \eqref{eq:bound_tilde_k_j}. 
\hfill $\blacksquare$

Next, we present the main theorem of the paper.

\textbf{Theorem 1.} Consider Algorithm 1 to solve $x$ for the proposed problem~\eqref{eq:overall_obj_fun}. Then a real $\rho_{\omega} \in (0, 1)$ and selections of $\beta \geq 0$ and $K(0)$ exist, such that for each iteration $t \geq 0$,
\begin{align}
\| g(t+1) \| \leq \left( \frac{\beta}{\lambda_r + \beta} + \lambda_1 \rho_{\omega}^{t+1} \| K(0) - K_{\beta} \|_F \right) \| g(t) \|,
\label{eq:theorem_condition}
\end{align}
where $\lambda_1$ and $\lambda_r$ are the largest and smallest nonzero eigenvalue of $\mathbf{A}^T \mathbf{A}$, respectively. Moreover, it leads to $\lim_{t \to \infty} \| g(t) \| = 0$, indicating linear convergence of Algorithm 1.

\textit{Proof:}
Let the error be defined as $z(t) = x(t) - x^*$, then the updating equation for $z(t+1)$ is
\begin{align}
    z(t+1) & = x(t+1) - x^* \nonumber \\
    & = x(t) - K(t + 1) \textstyle\sum_{i=1}^N g_i(t) - x^* \nonumber \\
    & = \left( I - K(t + 1) \mathbf{A}^T \mathbf{A} \right) z(t).
\end{align}

For an arbitrary solution $x^* \in X^*$, $\mathbf{A}^T (\mathbf{A} x^* - \mathbf{b}) = \mathbf{0}_d$, hence, $g(t) = \mathbf{A}^T \mathbf{A} z(t)$. Let $\mathcal{N}(\mathbf{A}^T \mathbf{A})$ and $\mathcal{N}(\mathbf{A}^T \mathbf{A})^\perp$ denote the nullspace of $\mathbf{A}^T \mathbf{A}$ and its orthogonal vector space, respectively. By this definition, we can write $z(t) = z(t)^\perp + z(t)^{\mathcal{N}}$, where $z(t)^\perp \in \mathcal{N}(\mathbf{A}^T \mathbf{A})^\perp$ and $z(t)^{\mathcal{N}} \in \mathcal{N}(\mathbf{A}^T \mathbf{A})$ are two orthogonal vectors. Due to $\mathbf{A}^T \mathbf{A} z(t)^{\mathcal{N}} = \mathbf{0}_d$, we can rewrite $g(t)$ as $g(t) = \mathbf{A}^T \mathbf{A} z(t)^\perp$.

Since $\mathbf{A}^T \mathbf{A} \succeq 0$, in general $1 \leq r \leq d$, where $r$ is the rank of $\mathbf{A}^T \mathbf{A}$. Recall the eigenvalues $ \lambda_1 \geq \dots \geq \lambda_d \geq 0 $. Let $S = \text{Diag}([\lambda_1, \dots, \lambda_d]^T)$,
and the matrix $ V $ consist of the corresponding orthonormal eigenvectors $ [V_1, \dots, V_d] $. Hence, $ V^T V = I $, and $ \mathbf{A}^T \mathbf{A} = V S V^T $. 

If $ r < d $, then $\lambda_1 \geq \dots \geq \lambda_r > \lambda_{r+1} = \dots = \lambda_d = 0$, and 
$S = \text{Diag}([\lambda_1, \dots, \lambda_r, 0, \dots, 0)]^T$.
Let $ \text{span}\{ V_1, \dots, V_r \} $ denote the vector space spanned by the orthonormal eigenvectors $ V_1, \dots, V_r $, i.e.,
\begin{align}
    \text{span}\{ V_1, \dots, V_r \} = \left\{ \textstyle\sum_{i=1}^r u_i V_i : u_i \in \mathbb{R} \right\}.
\end{align}
As the eigenvectors $ V_1, \dots, V_d $ are orthogonal, we have $\mathcal{N}(\mathbf{A}^T \mathbf{A})^\perp = \text{span} \{ V_1, \dots, V_r \}$ and $\mathcal{N}(\mathbf{A}^T \mathbf{A}) = \text{span} \{ V_{r+1}, \dots, V_d \}$. Let $S^\perp = \text{Diag}([\textbf{1}_r^T, \textbf{0}_{d-r}^T]^T)$. If we define a projection matrix $Q = V S^\perp V^T$, then it can be concluded that $z(t)^\perp = Q z(t)$. Hence,
\begin{align}
& z(t + 1)^\perp = Q z(t + 1)  = Q \left( I - K(t + 1) \mathbf{A}^T \mathbf{A} \right) z(t) \nonumber \\
& = z(t)^\perp - Q K(t + 1) \mathbf{A}^T \mathbf{A} z(t) \nonumber \\
& = z(t)^\perp - Q K(t + 1) \mathbf{A}^T \mathbf{A} \left( z(t)^{\mathcal{N}} + z(t)^\perp \right).
\end{align}

Since $\mathbf{A}^T \mathbf{A} z(t)^{\mathcal{N}} = 0$,
\begin{align}
& z(t + 1)^\perp = z(t)^\perp - Q K(t + 1) \mathbf{A}^T \mathbf{A} z(t)^\perp \nonumber \\
= & z(t)^\perp - Q K(t + 1) g(t).
\end{align}

By left multiplying $\mathbf{A}^T \mathbf{A}$ on both sides, and $g(t) = \mathbf{A}^T \mathbf{A} z(t)^\perp$, we have
{\small
\begin{align}
& g(t + 1) = g(t) - \mathbf{A}^T \mathbf{A} Q K(t + 1) g(t) \nonumber \\
= & g(t) - \mathbf{A}^T \mathbf{A} Q \left( \widetilde{K}(t+1) + K_{\beta} \right) g(t) \nonumber \\
= & \left( I - \mathbf{A}^T \mathbf{A} Q K_{\beta} \right) g(t) - \mathbf{A}^T \mathbf{A} Q \widetilde{K}(t+1) g(t).
\end{align}
}

Using triangle inequality,
{\small
\begin{align}
\| g(t + 1) \| & \leq \| \left( I - \mathbf{A}^T \mathbf{A} Q K_{\beta} \right) g(t)\| \nonumber \\
& \quad + \| \mathbf{A}^T \mathbf{A} Q \widetilde{K}(t+1) g(t)\| .
\label{eq:g_ineq}
\end{align}
}

First, the upper bound for the first term of~\eqref{eq:g_ineq} is derived. Since $ \mathbf{A}^T \mathbf{A} = V S V^T $ and $ V^T V = I $, we have
\begin{align}
    K_{\beta} & = V \, \text{Diag}([\frac{1}{\lambda_1 + \beta}, \dots, \frac{1}{\lambda_d + \beta}]^T) \, V^T, \\
    \mathbf{A}^T \mathbf{A} Q & = V \, \text{Diag}([\lambda_1, \dots, \lambda_r, \textbf{0}_{d-r}^T]^T) \, V^T.
\end{align}
Hence,
{\small
\begin{align}
    & I - \mathbf{A}^T \mathbf{A} Q K_{\beta} \nonumber \\
    = & V \, \text{Diag}([1-\frac{1}{\lambda_1 + \beta}, \dots, 1-\frac{1}{\lambda_r + \beta}, \textbf{1}_{d-r}^T]^T) \, V^T.
\end{align}
}
Since $g(t) \in \mathcal{N}(\mathbf{A}^T \mathbf{A})^\perp$, and $\lambda_1 \geq \dots \geq \lambda_r > 0$, we have

\scalebox{0.9}{\parbox{1\linewidth}{
\begin{align}
& \| \left( I - \mathbf{A}^T \mathbf{A} Q K_{\beta} \right) g(t)\| \nonumber \\
\leq & \left\| V \, \text{Diag}([1-\frac{1}{\lambda_1 + \beta}, \dots, 1-\frac{1}{\lambda_r + \beta}, \textbf{1}_{d-r}^T]^T) \, V^T \right\| \| g(t) \| \nonumber \\
\leq & \max \left\{ \left| 1 - \frac{\lambda_1}{\lambda_1 + \beta} \right|, \dots, \left| 1 - \frac{\lambda_r}{\lambda_r + \beta} \right| \right\} \| g(t) \|  \nonumber \\
\leq & \frac{\beta}{\lambda_r + \beta} \| g(t) \|.
\end{align}
}}

Next, we calculate the upper bound for the second term of~\eqref{eq:g_ineq}. By Lemma 2, for each $j \in \{1, \dots, d\}$, $\left\| \widetilde{k}_j(t) \right\|^2 \leq \rho_{\omega}^{2t} \left\| \widetilde{k}_j(0) \right\|^2$. Since $\left\| \widetilde{K}(t) \right\|_F^2 = \sum_{j=1}^{d} \left\| \widetilde{k}_j(t) \right\|^2$ and $\| \widetilde{K}(t) \| \leq \| \widetilde{K}(t) \|_F$~\cite{Horn_Johnson_1985}, for all $t \geq 0$ 
we have 
\begin{align}
\resizebox{0.3\textwidth}{!}{$
\| \widetilde{K}(t) \| \leq \left\| \widetilde{K}(t) \right\|_F \leq \rho_{\omega}^t \left\| \widetilde{K}(0) \right\|_F.
$}
\end{align}

Also, we have $\| \mathbf{A}^T \mathbf{A} Q \| = \lambda_1$~\cite{Horn_Johnson_1985}, hence, the upper bound of the second term of \eqref{eq:g_ineq} can be written as,
\begin{align}
\| \mathbf{A}^T \mathbf{A} Q \widetilde{K}(t+1) g(t)\| \leq \lambda_1 \rho_{\omega}^{t+1} \| \widetilde{K}(0) \|_F \| g(t) \|.
\end{align}

Therefore, by combining the two terms, \eqref{eq:g_ineq} becomes
\begin{align}
\| g(t + 1) \| & \leq \left( \frac{\beta}{\lambda_r + \beta} + \lambda_1 \rho_{\omega}^{t+1} \| \widetilde{K}(0) \|_F \right) \| g(t) \|,
\end{align}
which refers to~\eqref{eq:theorem_condition}.

Finally, we show that $\lim_{t \to \infty} \| g(t) \| = 0$. Let $\alpha(t) = \frac{\beta}{\lambda_r + \beta} + \lambda_1 \rho_{\omega}^{t+1} \| \widetilde{K}(0) \|_F$. Since $\rho_{\omega} \in (0, 1)$, $\alpha(t)$ is strictly decreasing, and $\lim_{t \to \infty} \alpha(t) = \frac{\beta}{\lambda_r + \beta} < 1$. Hence, at the iteration $\tau$ when $\alpha(\tau) < 1$ for the first time, we have
\begin{align}
    \| g(\tau + 1) \| & \leq \alpha(\tau) \| g(\tau) \| < \alpha(0) \| g(\tau) \| \nonumber \\
    & < \alpha(0)^{\tau+1} \| g(0) \|. \nonumber
\end{align}
Then for all $t>\tau$, 
\begin{align}
    g(t+1) & \leq \textstyle \prod_{k=\tau+1}^{t} \alpha(k) \alpha(0)^{\tau+1} \| g(0) \| \nonumber \\
    & \leq \alpha(\tau+1)^{t+1} (\frac{\alpha(0)}{\alpha(\tau+1)})^{\tau+1} \| g(0) \|,
\end{align}
which indicates that $\lim_{t \to \infty} \| g(t) \| = 0$. \hfill $\blacksquare$

\revision{To summarize, Theorem 1 formalizes convergence of Algorithm 1 via an iterative scheme that advances classical SOR preconditioning.} The algorithm converges in a linear rate in a general case when $\mathbf{A}^T \mathbf{A}$ can be rank deficient. Next, we consider $\mathbf{A}^T \mathbf{A}$ being full rank, which means $x^*$ is the unique solution for problem of~\eqref{eq:overall_obj_fun}. In this case, Algorithm 1 with $\beta = 0$ converges \textit{superlinearly} to $x^*$, which is proved next. \textcolor{black}{We emphasize that $\mathbf{A}^T \mathbf{A}$ being full rank does not imply each $A^i$ being full column rank, since $x^T \mathbf{A}^T \mathbf{A} x = \sum_{i=1}^{N} x^T A_i^T A_i x = \sum_{i=1}^{N} ||A_i x||^2 > 0$ for any $x \ne \mathbf{0}_d$, meaning $\cap_{i=1}^{N} \mathcal{N}(A_i) = \{\mathbf{0}_d\}$ and not necessarily each $\mathcal{N}(A_i) =  \{\mathbf{0}_d\}$. Thus each local cost $f_i$ still remains non-strongly convex.}

\textbf{Corollary 2. (Superlinear Convergence)} If a unique solution $x^*$ exists for the linear LS problem~\eqref{eq:overall_obj_fun}, by choosing $\beta = 0$, Algorithm 1 converges \textit{superlinearly} to $x^*$. 

\textit{Proof:} When a unique solution $x^*$ exists, the matrix $\mathbf{A}^T \mathbf{A}$ is full-rank. Hence, with $\beta = 0$ we have 
\begin{align}
   & \| g(t + 1) \|  \leq \lambda_1 \rho_{\omega}^{t+1} \| {K}(0) - \mathbf{A}^T \mathbf{A} \|_F \| g(t) \| \nonumber \\
\Rightarrow \quad &  \frac{\| g(t + 1) \|}{\| g(t) \|}  \leq \lambda_1 \rho_{\omega}^{t+1} \| {K}(0) - \mathbf{A}^T \mathbf{A} \|_F \nonumber \\
\Rightarrow \quad &  \lim_{t \to \infty} \frac{\| g(t + 1) \|}{\| g(t) \|} = 0,
\end{align}
which indicates Algorithm 1 converges superlinearly to $x^*$. \hfill $\blacksquare$

\begin{table*}[ht!]
\centering
\caption{\revision{The number of iterations required by different algorithms to attain error $\varepsilon_{tol}$.}}
\renewcommand{\arraystretch}{1.2}
\setlength{\tabcolsep}{6pt}
\begin{tabular}{|l|c|c|c|c|c|c|c|c|c|}
    \hline
    \textbf{Dataset} & $\lambda_1/\lambda_r$ & $\varepsilon_{tol}$ & Dist GD & NAG & APC & \revision{Dist L-BFGS} & \revision{ADMM} & Dist IPG & Algo. 1 \\
    \hline
    \textbf{ash608} & 11.38 & $10^{-4}$ & $> 10^2$ & 29 & 20 & \revision{19} & \revision{99} & 11 & \textbf{6} \\ 
    \textbf{bcsstm07} & $5.8 \times 10^7$ & $10^{-4}$ & $> 10^5$ & $5.64 \times 10^4$ &  $4.85 \times 10^4$ & \revision{$2.64 \times 10^4$} & \revision{$1.04 \times 10^4$} & $1.19 \times 10^4$ & \textbf{404} \\
    \textbf{gr\_30\_30} & $3.79 \times 10^4$ & $10^{-4}$ & $> 10^5$ & $1.94 \times 10^3$ & $1.11 \times 10^3$ & \revision{448} & \revision{$3.28 \times 10^4$}  & 697 & \textbf{410}\\
    \hline
\end{tabular}
\label{tab:iterations}
\end{table*}

\section{Numerical Examples}
In this section, the proposed algorithm is validated using benchmark problems from the SuiteSparse Matrix Collection~\cite{SuiteSparseMatrix} (http://sparse.tamu.edu). 
\revision{
These problems are originated from real engineering problems, which are typically sparse and best solved efficiently.} 
In each problem, the collective data matrix $\mathbf{A}$ and collective observation vector $\mathbf{b}$ are provided to satisfy $\mathbf{A}x^* = \mathbf{b}$, where $x^*$ is the parameter vector with all entries equal to $\mathbf{1}$. It is worth noting that for some problems $x^* = \mathbf{1}$ is the unique solution, while others have multiple feasible solutions. We evaluate the convergence performance by the iteration counts to achieve a predetermined tolerance value, which is an estimation error defined as 
$\varepsilon_{tol} = {\| x(t) - x^* \|}$ for unique solution, or 
$\varepsilon_{tol} = {\| \mathbf{A} x(t) - \mathbf{b} \|}$ for multiple solutions. 
\revision{Although the examples have exact solutions (either unique or multiple), the proposed approach still converges to solve~\eqref{eq:overall_obj_fun} when no exact solution of $Ax=b$ exists.}

First, we test Algorithm 1 for three problems, namely ``\textit{ash608}'', ``\textit{bcsstm07}'', and ``\textit{gr\_30\_30}''. The results are compared with those of other methods, including distributed GD, NAG, APC,
\revision{distributed L-BFGS, ADMM}, and distributed {\em IPG}~\cite{chakrabarti2022iterative}. The collective data are distributed into $N = 10$ agents in all the tests. Optimal parameters for other methods can be found in~\cite{chakrabarti2022iterative},~\cite{lessard2016analysis} and~\cite{azizan2019distributed} for the best performance. Note that these parameters require certain calculations in other methods, but such calculation is not needed in Algorithm 1. For Algorithm 1, each entry of $x(0)$ and $K(0)$ is initialized as zero. The parameters are selected as $\beta = 0.01$ and $\omega = 1.2$.

\begin{figure}[ht!]
    \centering
    \includegraphics[width=\linewidth]{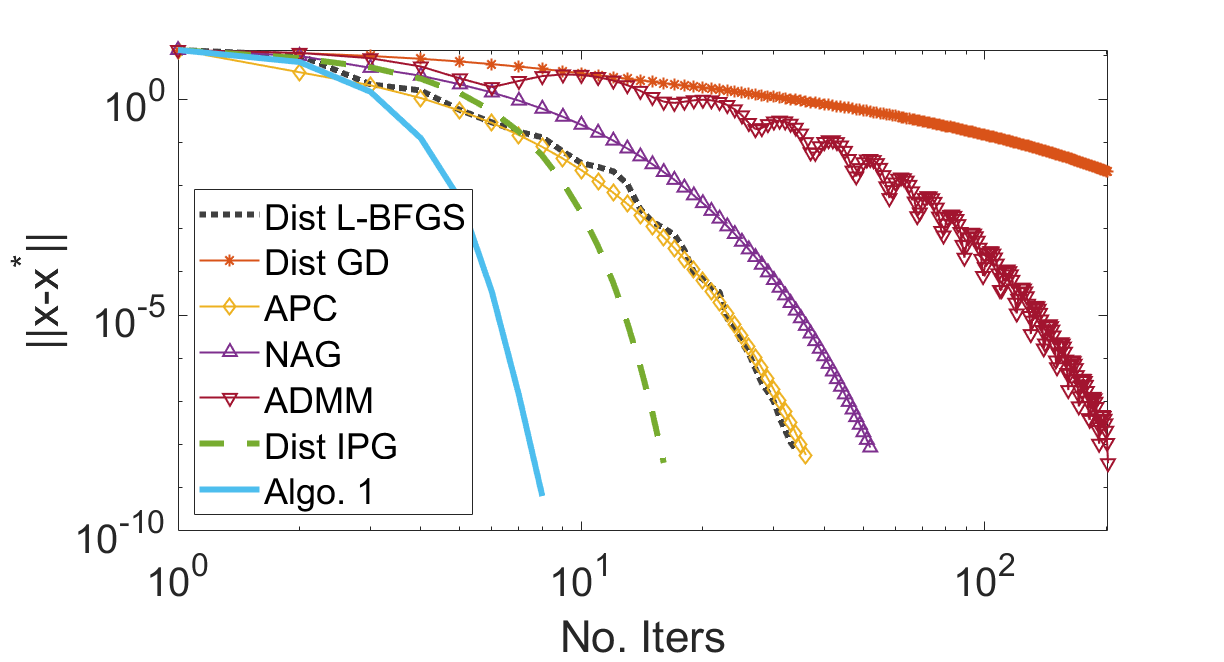}
    \caption{\revision{Temporal evolution of estimation errors in ``ash608''. (`Dist' means the distributed version of the corresponding algorithm.)}}
    \label{fig:ash608}
\end{figure}

The numbers of iterations to attain the required estimation error by different algorithms are listed in Table~\ref{tab:iterations}. It can be observed that the proposed Algorithm 1 outperforms existing methods, as well as our previous approach in~\cite{chakrabarti2022iterative}. Specifically, the convergence results are visualized in Fig.~\ref{fig:ash608} for the problem ``\textit{ash608}'', where the curves of the residue norms $\|x - x^*\|$ are plotted for all algorithms. 
\revision{In these problems ADMM requires a large number of iterations to converge, because the consensus mechanism must coordinate all agents, which results in a relatively slow convergence.}

Next, we test benchmark problems with large condition numbers. {\em IPG-SOR} iteratively updates a preconditioner matrix to accelerate the convergence, hence it is conjectured that the proposed method works effectively for problems with large condition numbers. $\varepsilon_{tol} = 10^{-4}$ is selected. The iteration count results are listed in Table~\ref{tab:large_cond}. 
\revision{As is observed, the counts for convergence are significantly reduced compared with distributed L-BFGS and ADMM.}

\begin{table}
\centering
\caption{\revision{The number of iterations required by different algorithms to attain relative estimation error $\varepsilon_{tol}$.}}
\setlength{\tabcolsep}{1pt}
\begin{tabular}{|l|c|c|c|c|c|c|c|c|}
    \hline
    \textbf{Dataset} & Cond. No. & \revision{Dist L-BFGS} & \revision{ADMM} & Dist IPG & Algo. 1\\
    \hline
    \textbf{abb313}     & $1.83 \times 10^{18}$ & \revision{275} & \revision{4327}           & $ 226 $   & $\textbf{52}$ \\
    \textbf{can\_229}   & $2.12 \times 10^{17}$ & \revision{414} & \revision{$1.63 \times 10^{4}$}  & $ 357 $   & $\textbf{76}$\\
    \textbf{Maragal\_2} & $2.86 \times 10^{47}$ & \revision{482} & \revision{889} & $ 321 $   & $\textbf{116}$\\
    \hline
\end{tabular}
\label{tab:large_cond}
\end{table}

\begin{table}[ht!]
\centering
\caption{\revision{Problem data size \& comparison of average computation time (milliseconds) per iteration per agent of different algorithms.}}
\renewcommand{\arraystretch}{1.2}
\setlength{\tabcolsep}{2pt}
\begin{tabular}{|l|c|c|c|c|c|c|c|c|}
    \hline
    \textbf{Dataset} & $\sum n_i$ & $d$ & Dist L-BFGS & ADMM & Dist IPG & Algo. 1 \\
    \hline
    \textbf{ash608}     & 608 & 188 & 0.0891 & 0.0392  & 1.0663  & 1.0262 \\
    \textbf{bcsstm07}   & 420 & 420 & 0.0143 & 0.1802  & 5.9159  & 6.0182 \\
    \textbf{gr\_30\_30} & 900 & 900 & 0.0241 & 1.2232  & 22.7211 & 28.1605 \\
    \textbf{abb313}     & 313 & 176 & 0.0075 & 0.0262  & 0.7317  & 0.9687 \\
    \textbf{can\_229}   & 229 & 229 & 0.0158 & 0.0447  & 1.4692  & 1.4438 \\
    \textbf{Maragal\_2} & 555 & 350 & 0.0262 & 0.0926  & 3.6447  & 3.5825 \\
    \hline
\end{tabular}
\label{tab:runtime}
\end{table}

\revision{The average computation time per iteration of different algorithms for the test problems are listed in Table~\ref{tab:runtime}. Experiments are performed on a laptop with Intel i9-13900HX CPU. Note that all computations are executed serially in current scripts without any parallelization implementation. It can be seen that the computational time of distributed {\em IPG} and {\em IPG-SOR} are comparable, but much longer than either distributed L-BFGS or ADMM. Hence, we conclude that when communication costs are negligible, distributed L-BFGS is the preferred choice due to its low computational complexity per iteration. However, for the scenarios where the communication overhead dominates the computational cost (e.g., each data transfer needs a long period), the reduction in iteration count can play a key role in improving the performance.}

\textcolor{black}{Finally, it is worth noting that although the multi-agent linear LS problems tested in this section possess exact solutions, our proposed approach can also handle linear LS without exact solutions. In such problems, the optimal objective function value $F(x^*)$ will not decrease to $0$, but the solution still satisfies that the total gradient $g^*=\mathbf{0}_d$, thus the proof holds valid.}

\section{Conclusion}
In this paper, we propose a novel accelerated distributed optimization algorithm to solve multi-agent linear LS problems. The algorithm is built upon the iterative preconditioning framework, utilizing the successive over-relaxation technique to improve the convergence performance. Theoretical guarantee for convergence is provided. Using sparse matrix benchmark problems, we demonstrate and analyze the performance advantages of the proposed approach.



\bibliographystyle{IEEEtran}
\bibliography{dist_sor_ipg}

\end{document}